\documentclass[journal,twoside]{IEEEtran}
\usepackage{algorithmic}
\usepackage{algorithm}
\usepackage{mathrsfs,amsthm}
\IEEEoverridecommandlockouts
\makeatletter
\def\ps@headings{%
\def\@oddhead{\mbox{}\scriptsize\rightmark \hfil \thepage}%
\def\@evenhead{\scriptsize\thepage \hfil \leftmark\mbox{}}%
\def\@oddfoot{}%
\def\@evenfoot{}}
\makeatother 
\usepackage{graphicx,amssymb,stfloats,subfigure,multicol}
\usepackage{epsfig,epsf,psfrag,amssymb,amsfonts,latexsym,graphicx,mathrsfs,subfigure}
\usepackage{dsfont}
\usepackage{chngpage}
\usepackage[dvips]{color}
\usepackage{textcomp}
\usepackage{hyperref}
\usepackage[hyphenbreaks]{breakurl}
\usepackage[normalem]{ulem}
\usepackage{cleveref}
\usepackage{array,multirow,graphicx}
\usepackage{caption}

 \def\old#1{}    

\begin{document}

\title{Distributed Estimation of Oscillations in Power Systems: an Extended Kalman Filtering Approach}
\author{\Large Zhe Yu$^\dagger$, Di Shi$^\dagger$, Zhiwei Wang$^\dagger$, Qibing Zhang$^\ddagger$, Junhui Huang$^\ddagger$, and Sen Pan${^\dagger}{^\dagger}$

\thanks{\scriptsize
Z. Yu$^\dagger$, D. Shi$^\dagger$ and Z. Wang$^\dagger$ are with GEIRI North America, San Jose, CA 95134, USA. Q. Zhang$^\ddagger$ and J. Huang$^\ddagger$ are with State Grid Jiangsu Electric Power Company, 215 Shanghai Road, Gulou, Nanjing, 210024, China. Sen Pan${^\dagger}{^\dagger}$ is with GEIRI, Nanjing, 210024, China. Email:
{\tt\{zhe.yu,di.shi,zhiwei.wang\}@geirina.net, zqb\_312@163.com, xshhjh@126.com,  pansen@geiri.sgcc.com.cn}.  This work is funded by SGCC Science and Technology Program under contract No. 5455HJ160007.}}

\maketitle

\begin{abstract}
Online estimation of electromechanical oscillation parameters provides essential information to prevent system instability and blackout and helps to identify event categories and locations. We formulate the problem as a state space model and employ the extended Kalman filter to estimate oscillation frequencies and damping factors directly based on data from phasor measurement units. Due to considerations of communication burdens and privacy concerns, a fully distributed algorithm is proposed using diffusion extended Kalman filter. The effectiveness of proposed algorithms is confirmed by both simulated and real data collected during events in State Grid Jiangsu Electric Power Company.
\end{abstract}


\section{Introduction}

\IEEEPARstart{E}{lectromechanical} oscillations are observed in interconnected power systems after large disturbances. Poorly damped oscillations reduce margins of power systems and could cause system instability or blackout. Wide-area measurement system (WAMS) technology using phasor measurement units (PMUs) makes it possible to observe the phenomenon of oscillations and estimate parameters such as frequencies, damping factors and magnitudes. These parameters contain vital information about modes of the power system and help operators to identify event categories and locations.

Compared to conventional supervisory control and data acquisition (SCADA) systems, PMUs have higher sampling rate and are able to measure phase angles, which attract great interests and investments in the past decade. Department of Energy (DoE) has spent over \$328 million in aggregate on synchrophasor technology and related communications networks \cite{NASPI:14TR}. By 2015, there were over 1,700 PMUs on the North American power grid, covering the entire U.S. high-voltage transmission network \cite{DoE:16TR}. In Jiangsu province of China, more than 160 PMUs have been deployed, covering all 500kV and majority of 220kV substations. 

On one hand, abundant time stamped data from PMUs allow system operators to monitor the system in detail and better understand its dynamics.
On the other hand, high sampling rate of PMU creates huge burden to the communication infrastructure under the existing centralized control mechanism, in which all data are uploaded to and processed in the control center. This mechanism limits many potential applications of WAMS which require PMU to have higher reporting rates. In addition, a centralized mechanism is vulnerable to cyber attacks. An improved framework is needed which can help better fulfill PMU's capabilities given limited network transit and computation capacity.

\subsection{Summary of Results}
We develop a fully distributed algorithm to monitor system oscillations. In this framework, information is exchanged locally and computation is distributed to each PMUs or PDC.

We first extend the formulation from \cite{Yazdanian&Etal:2015TPS} to a multi-measurements framework. A nonlinear system is formulated, whose state includes frequencies, damping factors and magnitudes of various oscillation modes from each PMU. The frequencies and damping factors are assumed consistent across the entire system while the amplitudes and phasors of each PMU may be different. The observation of the nonlinear system is measurements from PMUs plus white noises.

Then we proposes a centralized EKF which can directly estimate oscillation frequencies and damping factors of multiple modes with online implementation. The proposed method requires much lower computational resource compared to conventional methods because of the recursive nature of EKF.

Furthermore, a fully distributed EKF framework is developed in which no central coordinator is needed, and PMUs communicate information only with neighbour(s). The EKF computation is carried out at each PMU based on local information to estimate oscillation parameters. A consensus is achieved by a diffusion process.

\subsection{Related Work and Organization}
There is expanding literature on oscillation detection and estimation, most of which is focusing on centralized mechanisms. Some of the well known methods include matrix-pencil method (MP) \cite{Sarkar&Pereira:1995APM,Bounou&Lefebvre&Malhame:1992TPS}, eigenvalue realization algorithm (ERA) \cite{Juang&Pappa:1985JGCD,Peterson:1995JGCD}, Hankel total least-squares (HTLS) \cite{Sanchez&Chow:1997TPS}, Hilbert-Huang transform (HHT) \cite{Ruiz&Etal:2005PSCC}, Prony methods \cite{Hauer&Demeure&Scharf:1990TPS,Trudnowski&Johnson&Hauer:1999TPS}, and extended Kalman filter (EKF) \cite{Yazdanian&Etal:2015TPS, Peng&Nair:2012TPS, Jiang&Matei&Baras:2010WSCS}. Most of these methods are not designed for online implementation and do not scale up well. An attempt of applying online estimators of Prony method can be found in \cite{Zhou&Etal:2010PESGM}.

Due to the explosion of data volume, distributed computing and data processing algorithms started to obtain increasing attention in recent years \cite{Kar&Hug:2012PESGM,Zhang&Chow:2012TPS}, while few works has been reported on distributed oscillation monitoring. Authors in \cite{Nabavi&Zhang&Chakrabortty:2015TSG} model the Prony method as a consensus optimization and applied alternating direction methods of multipliers (ADMM) and exam  the algorithm's performance with different communication environments. However, the proposed ADMM method still requires a central coordinator to implement.

The work closest to this paper is \cite{Jiang&Matei&Baras:2010WSCS}, in which the authors presented a consensus extended Kalman filter to indirectly estimate oscillation modes. In this work, however, oscillation frequencies and damping factors are directly estimated. Furthermore, measurement diffusion and state reduction are applied to enhance the performance of EKF.

The remainder of the paper is organized as follows. In section \ref{sec:II}, we formulate a nonlinear state space model, whose states include oscillation frequencies, damping factors and magnitudes. A centralized extended Kalman filter is applied in section \ref{sec:III}. Considering the data volume and comminution burden, we present a fully distributed EKF framework in section \ref{sec:IV}. Numerical results based on simulated and real data in section \ref{sec:V} confirm the desirable performance of proposed algorithms, and section \ref{sec:VI} concludes the paper.

\section{Problem Formulation} \label{sec:II}

As discussed in \cite{Hauer&Demeure&Scharf:1990TPS} and \cite{Korba:2007GTD}, electromechanical oscillations in power systems can be represented as a sum of some exponentially damped sinusoids. In a discrete framework, a system measurement $y[k]=[y_1[k],\cdots,y_M[k]]^T$ can be expressed as follows.
\begin{equation}\label{eqn:signal}
\begin{array}{l}
y_m[k]=\sum_{l=1}^LA_{l,m}\exp(-\frac{\sigma_lk}{f_s})\cos(\frac{\omega_lk}{f_s}+\phi_{l,m})+\varepsilon_{m}[k],
\end{array}
\end{equation}
where $y_m[k]$ is the measurement of the $m$th PMU at the $k$th time instant\footnote{In practice, each PMU may contain multiple channels. For notation simplicity, we assume each PMU contains single channel measurements. Results can be easily extended to the case of multi-channel measurements.}, $T$ the transpose operator, $M$ the number of PMUs, ${A_{l,m}\in \mathbb{R}}$ the amplitude, $L$ the number of oscillation modes, $\sigma_l$ the damping factor, $\omega_l$ the frequency, $\phi_{l,m}$ the phase angle, $f_s$ the sampling rate, and $\varepsilon_m[k]$ the measurement error. The measurement noise is assumed to be a white Gaussian noise with zero mean and a diagonal covariance matrix ${R_k=diag(R_{1,k},\dots,R_{M,k})}$. Measurements from different PMUs may have various amplitudes and phase angles, while frequencies and damping factors are assumed to be consistent across the system.

Inspired by \cite{Yazdanian&Etal:2015TPS}, we formulate a  nonlinear system whose states contain frequencies and damping factors of the oscillation modes. Consider a sinusoid signal as follows.
\[
\begin{array}{rl}
s_{l,m}[k]&\triangleq A_{l,m}\exp(-\frac{\sigma_lk}{f_s})\cos(\frac{\omega_lk}{f_s}+\phi_{l,m})\\
&=\exp(-\frac{\sigma_lk}{f_s})A_{l,m}[\cos(\frac{\omega_lk}{f_s})\cos(\phi_{l,m})\\
&\mathrel{\phantom{=}}-\sin(\frac{\omega_lk}{f_s})\sin(\phi_{l,m})]\\
&=\exp(-\frac{\sigma_lk}{f_s})[B^c_{l,m}\cos(\frac{\omega_lk}{f_s})+B^{s}_{l,m}\sin(\frac{\omega_lk}{f_s})],
\end{array}
\]
where ${B^c_{l,m}\triangleq A_{l,m}\cos(\phi_{l,m})}$ and ${B^s_{l,m}\triangleq -A_{l,m}\sin(\phi_{l,m})}$.
Consider the evolution of the sinusoid signal as follows.
\[\hspace{-1em}
\begin{array}{rl}
s_{l,m}[k+1]&=\exp(-\frac{\sigma_l(k+1)}{f_s})B^c_{l,m}\cos(\frac{\omega_l(k+1)}{f_s})\\
&\mathrel{\phantom{=}}+\exp(-\frac{\sigma_l(k+1)}{f_s})B^{s}_{l,m}\sin(\frac{\omega_l(k+1)}{f_s})\\
&=[B^c_{l,m}(\cos(\frac{\omega_lk}{f_s})\cos(\frac{\omega_l}{f_s})-\sin(\frac{\omega_lk}{f_s})\sin(\frac{\omega_l}{f_s}))\\
&\mathrel{\phantom{=}}+B^s_{l,m}(\sin(\frac{\omega_lk}{f_s})\cos(\frac{\omega_l}{f_s})+\cos(\frac{\omega_lk}{f_s})\sin(\frac{\omega_l}{f_s}))]\\
&\mathrel{\phantom{=}}\times\exp(-\frac{\sigma_l}{f_s})\exp(-\frac{\sigma_lk}{f_s}).
\end{array}
\]

Define system states as signal magnitudes, frequencies, and damping factors as follows.
\[
\begin{array}{cl}
x_{l,m}[k]&=\left[
\begin{array}{c}
x^c_{l,m}[k]\\
x^s_{l,m}[k]
\end{array}\right]\\
&=\left[
\begin{array}{c}
B^c_{l,m}\exp(-\sigma_lk/f_s)\cos(\omega_lk/f_s)\\
B^s_{l,m}\exp(-\sigma_lk/f_s)\sin(\omega_lk/f_s)
\end{array}\right],\\
\omega_l[k]&=\omega_l,\\
\sigma_l[k]&=\sigma_l.
\end{array}
\]
The state transition is presented as follows.
\begin{equation}\label{eqn:stateTransition}
\begin{array}{cl}
x^c_{l,m}[k+1]&=\exp(-\frac{\sigma_l[k]}{f_s})x^c_{l,m}[k]\cos(\frac{\omega_l[k]}{f_s})\\
&\mathrel{\phantom{=}}-\exp(-\frac{\sigma_l[k]}{f_s})x^s_{l,m}[k]\sin(\frac{\omega_l[k]}{f_s})\\
&\mathrel{\phantom{=}}+\epsilon_{l,m}^c[k],\\
x^s_{l,m}[k+1]&=\exp(-\frac{\sigma_l[k]}{f_s})x^c_{l,m}[k]\sin(\frac{\omega_l[k]}{f_s})\\
&\mathrel{\phantom{=}}+\exp(-\frac{\sigma_l[k]}{f_s})x^s_{l,m}[k]\cos(\frac{\omega_l[k]}{f_s})\\
&\mathrel{\phantom{=}}+\epsilon_{l,m}^s[k],\\
\omega_l[k+1]&=\omega_l[k]+\epsilon_{l}^\omega[k],\\
\sigma_l[k+1]&=\sigma_l[k]+\epsilon_{l}^\sigma[k],\\
\end{array}
\end{equation}
where $\epsilon$ is the system noise. The measurement equation (\ref{eqn:signal}) can be written as
\[
y_m[k]=\sum_{l=1}^L(x^c_{l,m}[k]+x^s_{l,m}[k])+\varepsilon_m[k].
\]

Denote the amplitude vector of modes measured by the $m$th PMU by ${a_m[k]=\big[x_{1,m}[k], \cdots, x_{L,m}[k]\big]^T}$. Define the state of the system as ${x[k]=[a_1[k],\cdots,a_M[k],\omega_1[k],\sigma_1[k],\cdots, \omega_M[k],\sigma_M[k]]^T}$, which has a dimension of $(2ML+2L)$-by-1 and is comprised of two parts. The first part, $a_1[k],\cdots,a_M[k]$,  is the magnitude of oscillation modes measured  at different PMUs. The second part, $\omega_l$ and $\sigma_l$, is the frequency and damping factor of each mode and this part is the consensus across different PMUs.

We can write the transition in a general form as follows.
\[
x[k+1]=f(x[k])+\epsilon[k]
\]
where the transition function $f(\cdot)$ is nonlinear and can be derived from equation (\ref{eqn:stateTransition}). We assume that $\epsilon[k]$ is a white Gaussian noise with zero mean and covariance matrix $Q_k$.

Given the system states and equation (\ref{eqn:signal}), we obtain the observation function as follows.
\[
y[k]=Hx[k]+\varepsilon[k]
\]
where

\[H=\left[
\begin{array}{c}
H_1\\
H_2\\
\vdots\\
H_M
\end{array}
\right]
=\left[
\begin{array}{ccccc}
1~1 \cdots&0~0 \cdots & \cdots & 0~0 \cdots& 0~0 \cdots\\
0~0 \cdots&1~1 \cdots & \cdots & 0~0 \cdots& 0~0 \cdots\\
\vdots& \vdots& \vdots & \vdots&\vdots\\
0~0 \cdots&0~0 \cdots & \cdots & 1~1 \cdots& 0~0 \cdots\\
\end{array}
\right].
\]
Here ${H_i\in\mathbb{R}^{1\times (2ML+2L)}}$ is the $i$th row of the observation matrix $H$. The ${(2(i-1)L+1)}$th to ${2iL}$th elements in $H_i$ are $1$s and others are zeros. Thus, the constructed system is summarized as follows.

\begin{equation}\label{eqn:systemEquation}
\begin{array}{l}
x[k+1]=f(x[k])+\varepsilon[k]\\
y[k]=Hx[k]+\epsilon[k]
\end{array}
\end{equation}

\section{Centralized Extended Kalman Filter}\label{sec:III}
\begin{figure}[ht]
\centering
\psfrag{x}{{\tiny $\hat{x}[k]$}}
\psfrag{y1}{{\tiny$y_1[k]$}}
\psfrag{y2}{{\tiny$y_2[k]$}}
\psfrag{y3}{{\tiny$y_3[k]$}}
\psfrag{y4}{{\tiny$y_4[k]$}}
\psfrag{PMU1}{{\tiny PMU 1}}
\psfrag{PMU2}{{\tiny PMU 2}}
\psfrag{PMU3}{{\tiny PMU 3}}
\psfrag{PMU4}{{\tiny PMU 4}}
  \includegraphics[width=.4\textwidth]{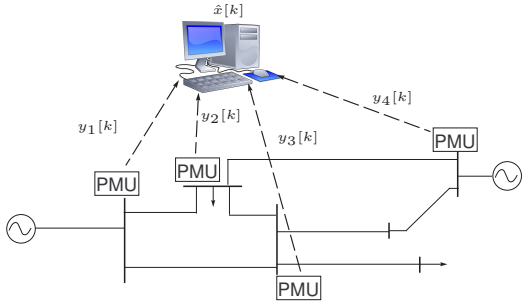}
  \caption{A centralized framework. PMUs send new measurements to a control center. The extended Kalman filter is carried out by the control center based on all data across the system.\normalsize}
\label{fig:centralized}
\end{figure}

In this section, a centralized framework of EKF is considered as shown in Figure \ref{fig:centralized}. At each time instant, new measurements are collected and sent from PMUs to a control center. The control center carries out a centralized extended Kalman filter to estimate the system state $x[k]$ based on all data across the system.

Given the system equations (\ref{eqn:systemEquation}), we apply an extended Kalman filter to estimate the system state. Kalman filter (KF) is a recursive algorithm to estimate the state of a linear dynamic system based on a series of noisy measurements. Based on the dynamic model, the KF predicts the priori state into the future and computes the difference between the predictions and the measurements. Then KF updates the posteriori estimation using the optimal Kalman gain and repeats the process. With white noises, Kalman
filter minimizes the mean squared estimation error.

When the dynamic system is nonlinear, the extended Kalman filter can be applied. Around the current estimated state, the EKF approximates the nonlinear system by a first order linearization and  applies the KF to the linearized system to find the optimal Kalman gain. The nonlinear system model and new measurements are used to calculate new state predictions. This process iterates and the state space model is re-linearized around updated state estimates.
\subsection{Centralized Extended Kalman Filter}
\begin{algorithm}[ht]
\caption{Centralized Extended Kalman Filter (CEKF)}
\label{alg:CEKF}
\begin{algorithmic}
\STATE $1$. Initialize $\hat{x}[0|-1]$ and $P[0|-1]$.
\STATE $2$. for $k=0:N-1$
~~~~\[
\begin{array}{l}
S=R_k+HP[k|k-1]H^T\\
K=P[k|k-1]H^TS^{-1}\\
\hat{x}[k|k]=\hat{x}[k|k-1]+K(y[k]-H\hat{x}[k|k-1])\\
P[k|k]=P[k|k-1]-KHP[k|k-1]\\
\hat{x}[k+1|k]=f(\hat{x}[k|k])\\
P[k+1|k]=F_kP[k|k]F_k^T+Q_k
\end{array}
\]
~~~end
\end{algorithmic}
\end{algorithm}
Let $\hat{x}[k|j]$ denote the minimum mean squared error estimate of $x[k]$ given measurements up to and including time $j$ and $P[k|j]$ the covariance matrix of the estimation error. Starting from the initial estimate $\hat{x}[0|-1]$ and $P[0|-1]$, the iteration of the extended Kalman filter for the system equation (\ref{eqn:systemEquation}) is summarized in Algorithm \ref{alg:CEKF}.

Here ${F_k=\frac{\partial f(x) }{\partial x}|_{x=\hat{x}[k|k]}}$ is the linearization of the system, and $N$ is the time length of measurements. The prediction process $f(\hat{x}[k|k])$ is stated as follows.
\[
f(\hat{x}[k|k])=\left[
\begin{array}{c}
\hat{a}_1[k+1|k]\\
\vdots\\
\hat{a}_{M}[k+1|k]\\
\hat{\omega}_1[k+1|k]\\
\hat{\sigma}_1[k+1|k]\\
\vdots\\
\hat{\omega}_L[k+1|k]\\
\hat{\sigma}_L[k+1|k]
\end{array}
\right]=\left[
\begin{array}{c}
\hat{x}_{1,1}[k+1|k]\\
\vdots\\
\hat{x}_{L,M}[k+1|k]\\
\hat{\omega}_1[k+1|k]\\
\hat{\sigma}_1[k+1|k]\\
\vdots\\
\hat{\omega}_L[k+1|k]\\
\hat{\sigma}_L[k+1|k]
\end{array}
\right]
\]
where
\[\small\begin{array}{l}
\hat{x}_{l,m}[k+1|k]\\
~~=\left[
\begin{array}{l}
\hat{x}^c_{l,m}[k+1|k]\\
\hat{x}^s_{l,m}[k+1|k]\\
\end{array}
\right]\\
~~=\left[
\begin{array}{l}
\exp(-\frac{\hat{\sigma}_l[k|k]}{f_s})[\hat{x}^c_{l,m}[k|k]\cos(\frac{\hat{\omega}_{i}[k|k]}{f_s})-\hat{x}^s_{l,m}[k|k]\sin(\frac{\hat{\omega}_{i}[k|k]}{f_s})]\\
\exp(-\frac{\hat{\sigma}_l[k|k]}{f_s})[\hat{x}^c_{l,m}[k|k]\sin(\frac{\hat{\omega}_{i}[k|k]}{f_s})+\hat{x}^s_{l,m}[k|k]\cos(\frac{\hat{\omega}_{i}[k|k]}{f_s})]
\end{array}
\right]\\
\hat{\omega}_{l}[k+1|k]=\hat{\omega}[k|k]\\
\hat{\sigma}_{l}[k+1|k]=\hat{\sigma}[k|k]
\end{array}
\]
\normalsize
\subsection{Initial Point and Coefficient Choice}
The accuracy and convergence of EKF rely heavily on the choice of initial points. In the context of oscillation estimation, a Fast Fourier Transform (FFT) or other similar technology can be employed as a trigger and the result can be used as a choice of initial points. FFT can estimate the spectral of sinusoids with limited measurements and alarm the operator with potential oscillations if the spectral of some frequency differs from noises significantly. Theses results can be used as inputs to Algorithm \ref{alg:CEKF}, and EKF will estimate the accurate frequency and damping factors. Other approaches such as singular value decomposition (SVD) \cite{Peng&Nair:2012TPS} can be applied to increase the confidence of the initial values.

Another possible choice is to use a look-up table, which can be built according to system operators' knowledge of the system and its typical oscillation modes. These modes can serve as the initial estimates which are fed into the EKF algorithm.

The proposed EKF algorithm is a model-based method and its performance relies on the proper choice of coefficients. Tuning of the covariance matrix of noise, $Q_k$ and $R_k$, is the major approach to adjust the performance of EKF. A large $Q_k$ or a small $R_k$ usually causes fluctuation around the actual value, while a small $Q_k$ or a large $R_k$ normally results in poor tracking. In this work, the tuning of coefficients is based on heuristic.

\section{Distributed Extended Kalman Filter}\label{sec:IV}
In this section, we consider a distributed framework of EKF. Note that, the term ``PMU'' in this section is in a much broader sense. It refers to any agent that can collect measurement data, communicate with other ``PMUs'', and carry out algorithms. A PMU in this work can be an actual PMU device, a PDC, a super PDC, or a data center.

\begin{figure}[ht]
\centering
\psfrag{x}{{\scriptsize $\hat{x}[k]$}}
\psfrag{y1}{{\scriptsize$y_1[k]$}}
\psfrag{y2}{{\scriptsize$y_2[k]$}}
\psfrag{y3}{{\scriptsize$y_3[k]$}}
\psfrag{y4}{{\scriptsize$y_4[k]$}}
\psfrag{PMU1}{{\scriptsize PMU 1}}
\psfrag{PMU2}{{\scriptsize PMU 2}}
\psfrag{PMU3}{{\scriptsize PMU 3}}
\psfrag{PMU4}{{\scriptsize PMU 4}}
\psfrag{PMU5}{{\scriptsize PMU 5}}
  \includegraphics[width=.45\textwidth]{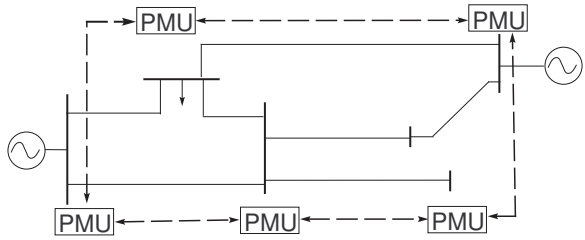}
  \caption{A fully distributed framework. PMUs communicate with their neighbour(s). EKF is carried out based on local information at each PMU. }
   \label{fig:distributed}
\end{figure}

 As shown in Figure \ref{fig:distributed}, no control center is needed in the fully distributed framework. At each time step, PMUs communicate with their neighbor(s), and EKF is carried out at each PMU. Under this framework, each PMU has its estimation, ${\hat{x}_m=[\hat{a}_{1,m},\cdots,\hat{a}_{M,m},\hat{\omega}_{1,m},\hat{\sigma}_{1,m},\cdots, \hat{\omega}_{M,m},\hat{\sigma}_{M,m}]^T}$, of the system state. The objective is to design algorithms to make the estimation converge to the actual value.

We formulate the topology of PMUs as a graph. Consider an undirected graph ${\mathcal{G}=\{\mathcal{N},\mathcal{E}\}}$, where ${\mathcal{N}=\{1,\dots,M\}}$ represents the PMU set and $\mathcal{E}$ the edge set. Each edge ${(i,j)\in\mathcal{E}}$ represents that PMU $i$ and PMU $j$ can communicate with each other. We define the set of nodes connected to a certain PMU $i$ as the neighbors of $i$, denoted by ${\mathcal{N}_i=\{j\in\mathcal{N}:(i,j)\in\mathcal{E}\}}$. A PMU is aways a neighbour of itself. The number of neighbors of PMU $m$ is referred to as degree, denoted by $|\mathcal{N}_m|$. Here we assume the graph is connected.

We extend a distributed Kalman filter framework proposed by \cite{Cattivelli&Sayed:TAC10}, referred to as diffusion Kalman filter. The diffusion Kalman filter attempts to approximate the global KF estimation by local information.
\begin{figure}[ht]
\centering
\subfigure[Exchange measurements]{\label{fig:diffusionEKF}
\psfrag{x}{{\scriptsize $\hat{x}[k]$}}
\psfrag{y1}{{\scriptsize$y_1[k]$}}
\psfrag{y2}{{\scriptsize$y_2[k]$}}
\psfrag{y3}{{\scriptsize$y_3[k]$}}
\psfrag{y4}{{\scriptsize$y_4[k]$}}
\psfrag{PMU1}{{\scriptsize PMU 1}}
\psfrag{PMU2}{{\scriptsize PMU 2}}
\psfrag{PMU4}{{\scriptsize PMU 4}}
\psfrag{PMU5}{{\scriptsize PMU 3}}
  \includegraphics[width=.22\textwidth]{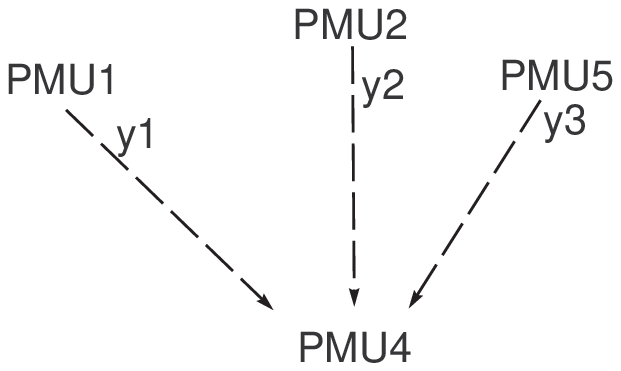}
  }
\subfigure[Exchange pre-estimates]{
\label{fig:exchangePreEstimate}
\psfrag{x}{{\scriptsize $\hat{x}[k]$}}
\psfrag{y1}{{\scriptsize$\varphi_4[k]$}}
\psfrag{y2}{{\scriptsize$\varphi_4[k]$}}
\psfrag{y3}{{\scriptsize$\varphi_4[k]$}}
\psfrag{PMU1}{{\scriptsize PMU 1}}
\psfrag{PMU2}{{\scriptsize PMU 2}}
\psfrag{PMU4}{{\scriptsize PMU 4}}
\psfrag{PMU5}{{\scriptsize PMU 3}}
  \includegraphics[width=.22\textwidth]{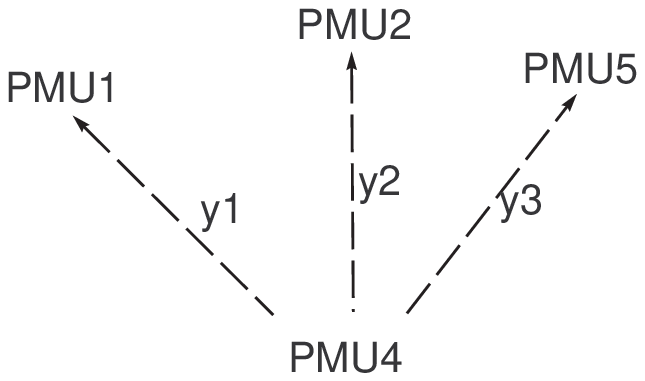}
  }
  \caption{Information exchanges at PMU 4.}

\end{figure}
As shown in Figure \ref{fig:diffusionEKF}, each PMU takes new measurements and collects new information from its neighbour(s). Based on this local information, each PMU carries out EKF to obtain a pre-estimation of the system state, $\varphi_m$. Then PMUs broadcasts its pre-estimation to its neighbour(s) and updates its estimate $\hat{x}_m$ by diffusion of all the pre-estimation collected from its neighbour(s). The diffusion EKF algorithm is described in Algorithm \ref{alg:DEKF}.

\begin{algorithm}
\caption{Diffusion EKF (DEKF)}
\label{alg:DEKF}
\begin{algorithmic}
\STATE $1$. Initialize $\hat{x}_m[0|-1]$ and $P_m[0|-1]$ for each PMU $m$.
\STATE $2$. for $k=0:N-1$\\
~~~~~~Each PMU sends measurements to its neighbour(s).\\
~~~~~~Incremental Update\\
~~~~~~~~~~for $m=1:M$\\
\[
\begin{array}{l}
\varphi_m=\hat{x}_{m}[k|k-1]\\
P_m=P_m[k|k-1]
\end{array}
\]
~~~~~~~~~~~~~for $j\in\mathcal{N}_m$\\
\[
\begin{array}{l}
S=R_{j,k}+H_jP_mH_j^T\\
K=P_mH_j^TS^{-1}\\
\varphi_m=\varphi_m+K(y_j[k]-H_j\varphi_m)\\
P_m=P_m-KH_jP_m
\end{array}
\]
~~~~~~~~~~~~~end\\
~~~~~~~~~~end\\
~~~~~~Each PMU sends pre-estimation to its neighbour(s).\\
~~~~~~Diffusion Update\\
~~~~~~~~~~for $m=1:M$\\
\[
\begin{array}{l}
\hat{x}_m[k|k]=\sum_{j\in\mathcal{N}_m}c_{m,j}\varphi_{j}\\
P_m[k|k]=P_m\\
\hat{x}_m[k+1|k]=f(\hat{x}_m[k|k])\\
P_m[k+1|k]=F_{m,k}P_m[k|k]F_{m,k}^T+Q_k
\end{array}
\]\\
~~~~~~~~~~end\\
~~~end
\end{algorithmic}
\end{algorithm}
Here ${F_{m,k}=\frac{\partial f(x) }{\partial x}|_{x=\hat{x}_m[k|k]}}$ , and $c_{m,j}$ is a diffusion factor and satisfies the following properties.
\begin{equation}
\label{eqn:diffusionFactors}
\begin{array}{l}
\sum_{j\in\mathcal{N}_m}c_{m,j}=1\\
c_{m,j}\ge0
\end{array}
\end{equation}
The diffusion of the pre-estimation is a weighted average. Note that, the diffusion update is not taken into account in the recursion of the matrices $P_m[k|k-1]$ and   $P_m[k|k]$, and they are no longer the covariance matrix of the state estimation.

\subsection{State Reduction}
In Algorithm \ref{alg:DEKF}, each PMU estimates the entire states of the system, including amplitudes of signals from PMUs that are not its neighbour(s). In the incremental update, for each PMU, only its neighbors' states are updated while the rest stay unchanged, which makes the process slow and the communication inefficient. Here we propose a state reduction framework of the DEKF to enhance the performance. Define the reduced estimate made by the $m$th PMU as ${\hat{x}^-_m=\big[\hat{a}_{m,j_1},\cdots,\hat{a}_{m,j_{|\mathcal{N}_m|}},\underline{\hat{x}}^-_{m}\big]^T}$, where ${\mathcal{N}_m=\{j_1,\cdots,j_{|\mathcal{N}_m|}\}}$ is the neighbor set of the $m$th PMU, $\hat{a}_{m,j}$ is the estimate of the $j$th PMU's amplitudes by the $m$th PMU, and $\underline{\hat{x}}^-_{m}$ is the estimate of the frequency and damping factor part. The reduced estimate is a $(2|\mathcal{N}_m|L+2L)$-by-1 dimension vector. Frequencies and damping factors of all modes and amplitudes of all neighboring PMUs are included while those of the non-neighboring ones are excluded.

In this case, the observation matrix $H$, estimation covariance matrix $P_m[k|i]$, process noise covariance $Q_k$, system functions $f(x)$ and Jacobian matrix $F_{m,k}$ are modified to $H_m^-$, $P_m^-[k|i]$, $Q^-_{m,k}$, $R^-_{m,k}$, $f^-_m(x)$ and $F_{m,k}^-$, accordingly. For the $m$th PMU, define the observation matrix ${H_m^-\in\mathbb{R}^{|\mathcal{N}_m|\times (2|\mathcal{N}_m|L+2L)}}$ as
\[\begin{array}{ll}
H_m^-&=\left[
\begin{array}{c}
H_{m,j_1}^-\\
H_{m,j_2}^-\\
\vdots\\
H_{m,j_{|\mathcal{N}_m|}}^-\\
\end{array}
\right]\\
&=\left[
\begin{array}{ccccc}
1~1\cdots& 0~0\cdots\cdots&\cdots&0~0\cdots& 0~0 \cdots\\
0~0\cdots& 1~1\cdots\cdots&\cdots&0~0\cdots& 0~0 \cdots\\
\vdots&\vdots&\vdots&\vdots &\vdots\\
0~0\cdots& 0~0\cdots\cdots&\cdots&1~1\cdots& 0~0 \cdots\\
\end{array}
\right],
\end{array}
\]
where $H^-_{m,j_i}$ is the $i$th row of $H_m^-$ with the ${(2(i-1)L+1)}$th to $2iL$th elements being $1$s and the rest being zeros. The formulas of other matrices are derived accordingly and omitted here.

The diffusion EKF under the reduced state framework is similar to Algorithm \ref{alg:DEKF}. Each PMU receives measurements from its neighbour(s) and estimates accordingly. After obtaining the pre-estimates, PMUs comminute this information and make a diffusion to update its estimate.  However, each PMU only maintains amplitude estimates of its neighbour(s) and the diffusion is carried out across neighbour(s) who estimates the same amplitudes.

Recall that the reduced estimation ${\hat{x}^-_m}$ by the $m$th PMU is comprised of the amplitudes of its neighbors and the frequency and damping factor part.
Denote $\varphi^a_{m,j}$ as the pre-estimate of $\hat{a}_{m,j}$, $\underline{\varphi}^-_{m}$ as the one of $\underline{\hat{x}}^-_{m}$, and $\varphi_m^-$ as the one of $\hat{x}^-_{m}$. The reduced state diffusion EKF is summarized in Algorithm \ref{alg:DEKF_reduce}.
\begin{algorithm}
\caption{Diffusion EKF with reduced state (DEKF-R)}
\label{alg:DEKF_reduce}
\begin{algorithmic}
\STATE $1$. Initialize $\hat{x}^-_m[0|-1]$ and $P^-_m[0|-1]$ for each PMU $m$.
\STATE $2$. for $k=0:N-1$\\
~~~~~~Each PMU sends  measurements to its neighbour(s).\\
~~~~~~Incremental Update\\
~~~~~~~~for $m=1:M$\\
\[
\begin{array}{l}
\varphi_m^-=\hat{x}^-_{m}[k|k-1],\\
P_m^-=P_m^-[k|k-1];
\end{array}
\]
~~~~~~~~~~~for $j\in\mathcal{N}_m$\\
\[
\begin{array}{l}
S=R_{j,k}+H_{m,j}^-P^-_m(H_{m,j}^-)^T,\\
K=P^-_m(H_{m,j}^-)^TS^{-1},\\
\varphi_m^-=\varphi_m^-+K(y_j[k]-H_{m,j}^-\varphi^-_m),\\
P_m^-=P_m^--KH_{m,j}^-P_m^-.
\end{array}
\]
~~~~~~~~~~~end\\
~~~~~~~~end\\
~~~~~~Each PMU sends pre-estimation to its neighbour(s).\\
~~~~~~Diffusion Update\\
~~~~~~~~for $m=1:M$\\
~~~~~~~~~~for $j\in\mathcal{N}_m$\\
\[
\hat{a}_{m,j}[k|k]=\sum_{i\in\mathcal{N}_m\cap\mathcal{N}_j}d_{m,j,i}\varphi^a_{i,j},\]
~~~~~~~~~~end\\
\[
\begin{array}{l}
\hat{\underline{x}}_m[k|k]=\sum_{j\in\mathcal{N}_m}c_{m,j}\underline{\varphi}^-_{j},\\
\hat{x}^-_m[k|k]=\\
~~~~\big[\hat{a}_{m,j_1}[k|k],\cdots,\hat{a}_{m,j_{|\mathcal{N}_m|}}[k|k],\underline{\hat{x}}^-_{m}[k|k]\big]^T,\\
P_m^-[k|k]=P_m^-,\\
\hat{x}^-_m[k+1|k]=f_m^-(\hat{x}^-_m[k|k]),\\
P^-_m[k+1|k]=\\
~~~~F^-_{m,k}P^-_m[k|k](F^-_{m,k})^T+Q^-_{m,k}.
\end{array}
\]
~~~~~~~~end\\
~~~end
\end{algorithmic}
\end{algorithm}

Here $d_{m,j,i}$ and $c_{m,j}$ are diffusion factors, where $c_{m,j}$ satisfies properties stated in equation (\ref{eqn:diffusionFactors}) and $d_{m,j,i}$ satisfies the following properties.
\[
\begin{array}{l}
\sum_{i\in\mathcal{N}_m\cap\mathcal{N}_j}d_{m,j,i}=1\\
d_{m,j,i}\ge0
\end{array}\]

\section{Numerical Results}\label{sec:V}
In this section, we present numerical results using both simulated and real PMU data collected from real-world system oscillation events. We first apply the proposed algorithms on a noisy ring down sinusoid signal and compare the accuracy of the proposed  algorithms with PRONY \cite{Trudnowski&Johnson&Hauer:1999TPS} and ADMM-PRONY \cite{Nabavi&Zhang&Chakrabortty:2015TSG}, a decentralized extension of PRONY. Then we test EKF and DEKF-R using a test case library  \cite{Sun&Etal:2016PESGM}. After that, real oscillation data from Jiangsu Electric Power Company in China are examined.

\subsection{Ring Down Sinusoids with Different Noise Levels}
In this case, the measurement is an exponentially damped sinusoid with a zero-mean white Gaussian noise stated as follows.
\[
y_m[k]=m\big[\exp(-\sigma k/f_s)\cos(\omega/f_sk+\phi_m)+\varepsilon_m[k]\big],
\]
where the frequency is ${\omega=4\pi}$ rad/s, the damping factor ${\sigma=0.0126}$, $\phi_m$ the phase angle, and $\varepsilon_m[k]$ the noise. The corresponding damping ratio is ${\zeta=\sigma/\sqrt{\sigma^2+\omega^2}=0.1\%}$. The phase angle $\phi_m$ is assumed to be uniformly distributed within $[-\pi/2,\pi/2]$. The total number of PMUs is set to be ${M=5}$. The amplitudes of PMUs are made different to model real-world signals from power systems. The sample rate $f_s$ is selected to be 30Hz and the length of the time window is set as $10$ seconds.

For the centralized EFK algorithm, the initial point is assumed to be uniformly distributed within a $[-70\%,130\%]$ range of the real value. The filter parameters are selected as $R_k=10^{-3}I$ where $I$ is an identity matrix with a proper dimension. $Q_k$ is a diagonal matrix whose first $2ML$ diagonal elements are set zeros and the last $2L$ diagonal elements are around $10^{-9}$.

For the distributed EKF algorithm, the communication topology is illustrated in Figure \ref{fig:distributed}, where the communication path of five PMUs forms a connected graph with each PMU communicating with its neighbor(s) only. The initial estimate at each PMU is assumed to be uniformly distributed within a $[-70\%,130\%]$ range of its corresponding true value. Diagonal elements in the measurement noise covariance matrix  $R_{m,k}$ are around $10^{-4}$ and the process noise covariance $Q_{m,k}$ is a diagonal matrix whose first $2|\mathcal{N}_m|L$ diagonal elements are zeros and last $2L$ diagonal elements are around $10^{-8}$.

The authors in \cite{Nabavi&Zhang&Chakrabortty:2015TSG} extended PRONY algorithm to a decentralized version, referred as ADMM-PRONY. Under this framework, a system coordinator is still needed who collects information across the entire system. Each PMU first carries out a linear estimation locally based on private measurements and sends the coordinator estimation results. The coordinator collects all estimates, averages them and broadcasts the diffusion back to PMUs. Given this global diffusion, each PMU solves a quadratic optimization to trade off the estimation accuracy and the tracking error of the diffusion, weighted by $\rho$, and sends the updated estimate to the coordinator. The iteration continues till the diffusion converges. Compared to the centralized PRONY, ADMM algorithms require  exchange of the estimates rather than the PMU measurements between PMUs and the coordinator, which relieves greatly the communication burden. However, the requirement of a system-level coordinator  makes it not fully distributed, and the system will be vulnerable and subject to a single point of failure at the coordinator. In this simulation, the weight of tracking error is set ${\rho=0.01}$ and the convergent tolerance is set as $0.01$.
\begin{table}
\caption{Single sinusoid  with different noise levels}
\centering
\begin{tabular}{|l|c|c|c|c|}
  \hline
 Error& Freq  & Damping &Freq  & Damping  \\

  \hline
  &\multicolumn{2}{c|}{SNR=50db}&\multicolumn{2}{c|}{SNR=40db}\\
  \hline
  \mbox{Mean (PRONY) } & .00\% & .48\% &.00\%& 1.51\%\\
  \mbox{Std (PRONY) } & .00\%& .38\%&.00\%& 1.17\%\\
  \hline
    \mbox{Mean (ADMM) } & .00\%& .84\%&.00\%&2.67\%\\
  \mbox{Std (ADMM) } & .00\% & .74\% &.00\% &2.03\%\\
  \hline
  \mbox{Mean (EKF) } & .00\% & .55\%   &.00\%&1.71\%\\
  \mbox{Std (EKF) } & .00\% & .42\%  &.00\%&1.24\% \\
  \hline
  \mbox{Mean (DEKF) } &.02\% &4.01\%&.07\%&7.29\% \\
  \mbox{Std (DEKF) } & .00\%& 2.73\%    &.02\%&5.63\%\\
  \hline
  \mbox{Mean (DEKF-R) } &.00\% &3.33\%&.05\%&6.9\% \\
  \mbox{Std (DEKF-R) } & .00\%& 2.41\%    &.03\%&4.05\%\\
  \hline
   &\multicolumn{2}{c|}{SNR=30db}&\multicolumn{2}{c|}{SNR=20db}\\
   \hline
  \mbox{Mean (PRONY) } & .00\% & 4.68\% &.02\%& 18.26\%\\
  \mbox{Std (PRONY) } & .00\%& 3.58\%&.01\%& 85.11\%\\
  \hline
    \mbox{Mean (ADMM) } & .00\%& 9.83\%&.03\%&64.85\%\\
  \mbox{Std (ADMM) } & .00\% & 116.04\% &.02\% &1148.29\%\\
  \hline
  \mbox{Mean (EKF) } & .00\% & 4.02\%   &.01\%&11.86\%\\
  \mbox{Std (EKF) } & .00\% & 2.98\%  &.01\%&9.04\% \\
  \hline
  \mbox{Mean (DEKF) } &.07\% &13.06\%&.07\%&35.68\% \\
  \mbox{Std (DEKF) } & .03\%& 10.05\%    &.05\%&29.33\%\\
  \hline
  \mbox{Mean (DEKF-R) } &.05\% &11.41\%&.05\%&35.07\% \\
  \mbox{Std (DEKF-R) } & .03\%& 8.59\%    &.03\%&26.28\%\\
  \hline
\end{tabular}
\label{table:singleSin}
\end{table}

One thousand Monte Carlo runs for each level of noise are carried out. Mean and standard deviations of estimation error are summarized in Table \ref{table:singleSin}. It is shown that, performance of centralized EKF is close to PRONY. The accuracy of distributed EKF with reduced state is comparable with the one of ADMM, while DEKF-R requires no coordinator and no global information exchange.
    With a low level of noise, both EKF and PRONY work well. DEKF-R gives slightly bigger error as compared to ADMM due to less information exchange, but the errors are all within the acceptable range. As the level of noise increases, EKF starts to outperform PRONY, both in the accuracy and stability of the estimates. ADMM is based on PRONY and also sensitive to noise and performs worse than DEKF when SNR is as large as 20db. Another observation to note is that DEKF-R dominates DEKF which indicates the effectiveness of state reduction.

\subsection{WECC Test Case Library}
Authors of \cite{Sun&Etal:2016PESGM} established a test case library for oscillation detection and forced oscillation source location in power systems.
 A reduced WECC 179-bus 29-machine system is simulated in TSAT with integration a step size around $0.004$s and sampled at rate $30$Hz. All generators are presented by a classical second-order differential model with damping parameter equals to $4$.

\begin{figure}[ht]
\centering
\subfigure[Rotor speed measurements]{
\label{fig:testCaseRotorSpd}
  \includegraphics[width=.5\textwidth]{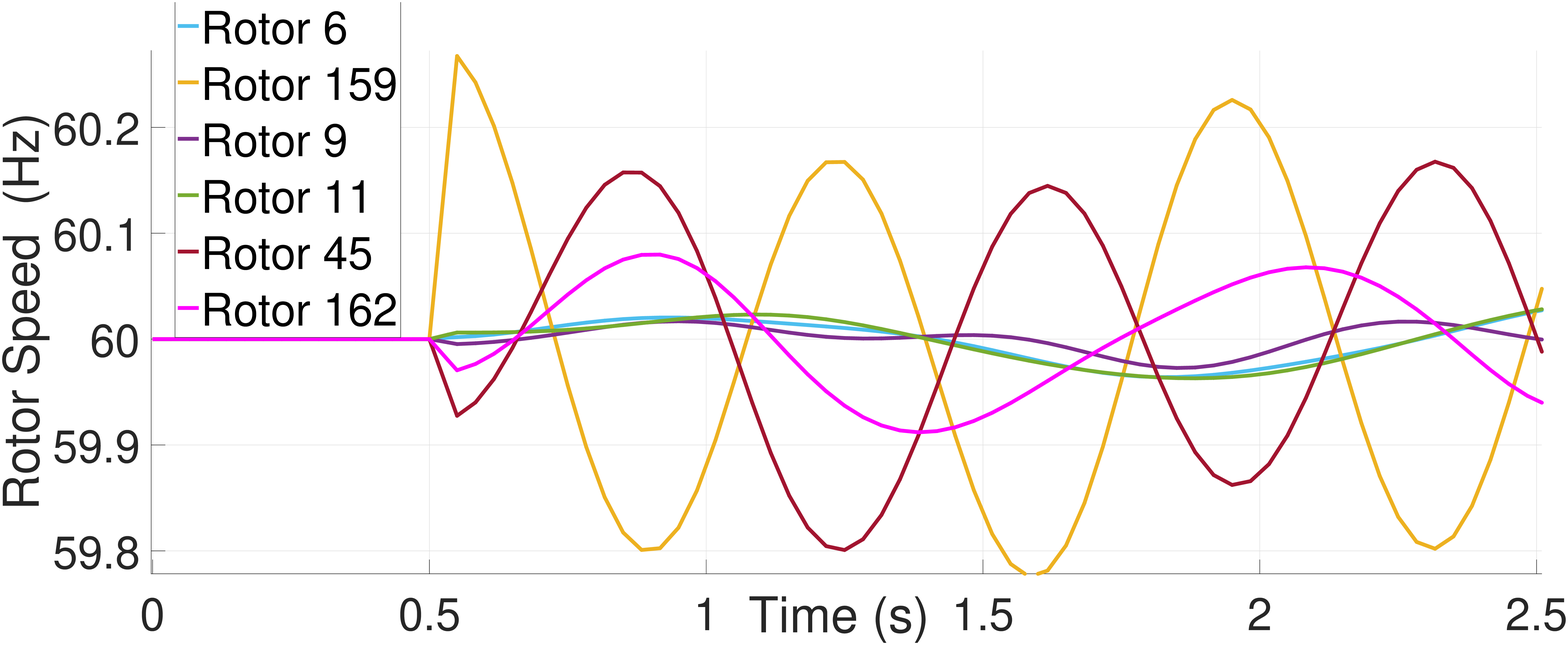}
  }
  \subfigure[FFT]{
\label{fig:testCase1FFT}
  \includegraphics[width=.5\textwidth]{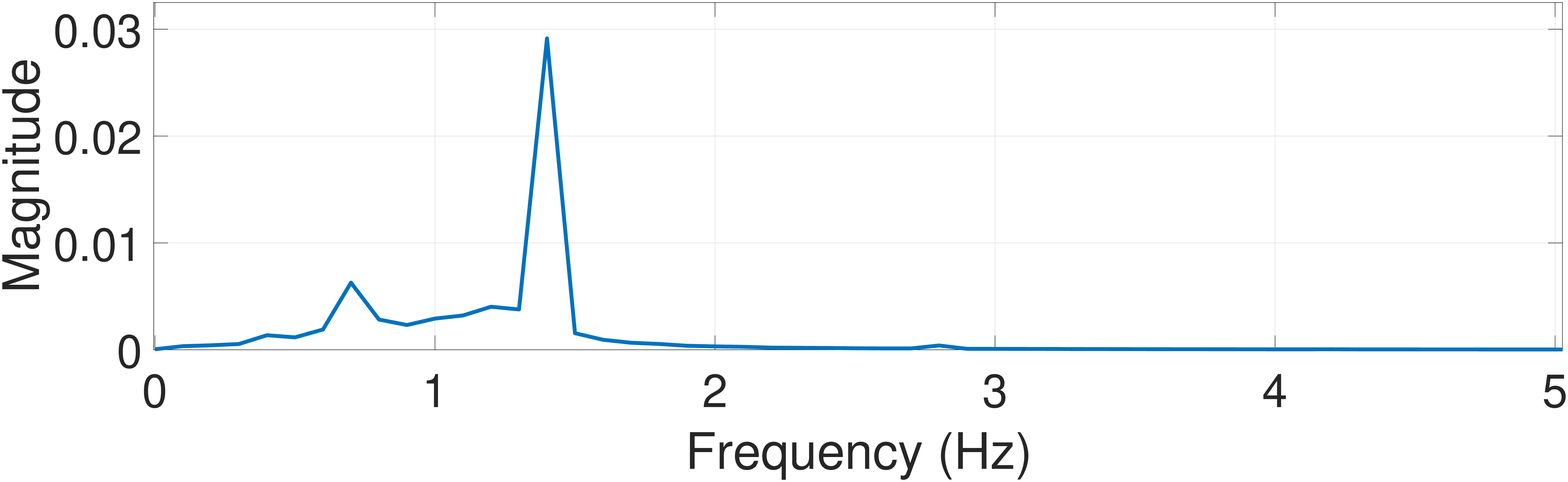}
  }
\subfigure[Measurements and fitted curve of rotor 159]{
\label{fig:testCaseEKF}
  \includegraphics[width=.5\textwidth]{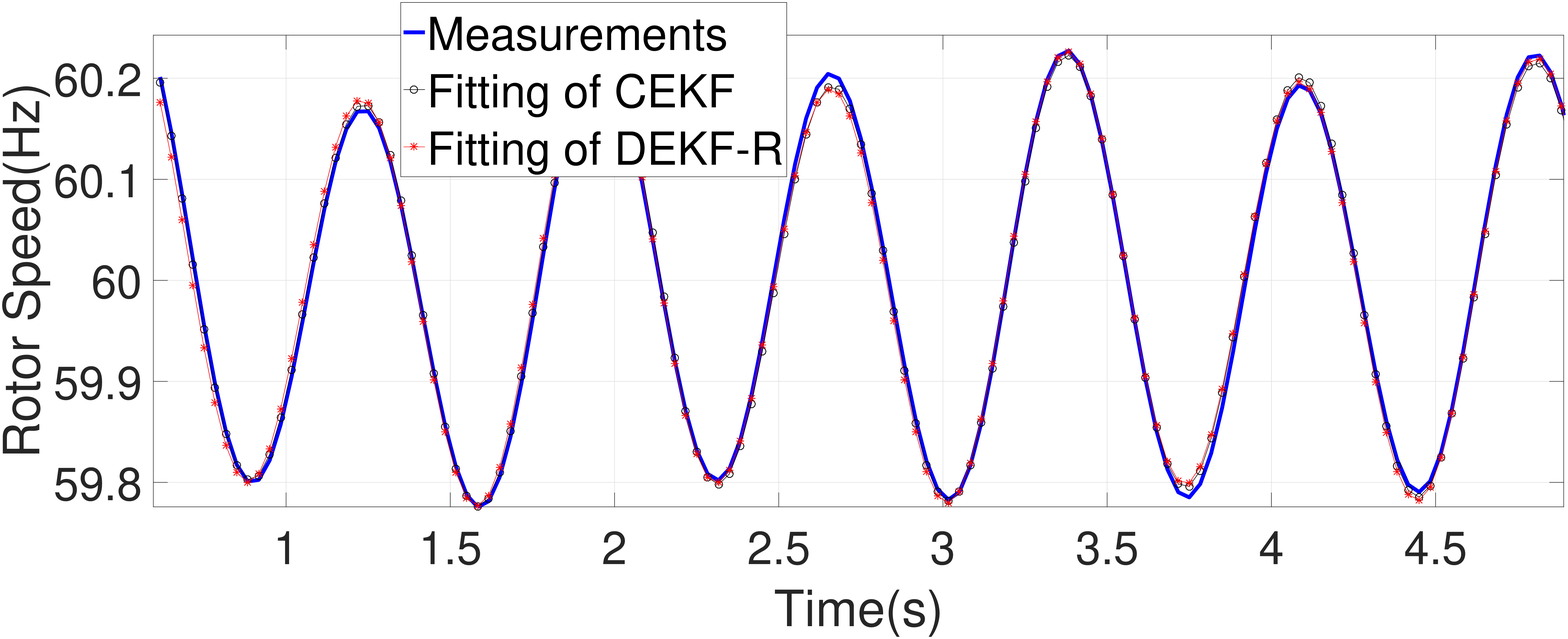}
  }
  \caption{WECC test case 1: rotor speed}
\label{fig:testCase1}
\end{figure}

In each test case, damping parameters of some generators are set such that they are poorly or negatively damped. Taking the first test case as an example, the damping factor of generator $45$ and $159$ are set to be $-2$ and $1$, respectively. At $0.5$ second, a  three-phase short circuit is added at bus $159$ and cleared by $0.55$ second to trigger oscillations in the system. As can be found in Figure \ref{fig:testCaseRotorSpd}, before $0.5$ second, the speed of rotors remains at $60$Hz. After the fault at $0.5$ second, all measurements begin to oscillate with different magnitudes and shift phase angles.

In this case, we first apply FFT to determine the number of modes and generate initial estimates for EKF methods. Feed FFT with speed data of rotor 30 from $0.7$ second to $10.7$ second, and the spectral shows there are two modes of frequency $0.7$Hz and $1.4$Hz, respectively. Set the initial estimates at ${\omega_1=0.7}$, ${\omega_2=1.4}$ and ${\sigma_1=\sigma_2=0}$ and run centralized EKF with parameters ${Q_k=10^{-9}I}$ and ${R_k=10^{-3}I}$. The estimation results are summarized in Table \ref{table:testCase1}. It can be found that EKF successfully identify the poorly damped frequency  at $1.4$Hz, and the difference between the measurements and fitted signal is small. The measurements and fitted speed curve of rotor $159$ are plotted in Figure \ref{fig:testCaseEKF} as an example.

\begin{table}[ht]
\caption{Test case 1: 2 oscillation modes}
\centering
\begin{tabular}{|l|c|c|c|c|}
  \hline
 Estimate& $\omega$ & $\sigma$ & $\omega$ & $\sigma$ \\
\hline
CEKF& $1.4016$ & $-0.0016$& $0.6927$ & $0.4715$\\
\hline
DEKF-R& $1.3999$ & $-0.0011$& $0.6892$ & $0.2761$\\
\hline
\end{tabular}

\label{table:testCase1}
\end{table}

Substitute the same initial points to the proposed Algorithm \ref{alg:DEKF_reduce} assuming that the communication of PMUs forms a circle like the one shown in Figure \ref{fig:distributed}, and each PMU communicates with its two neighbours. As shown in Figure \ref{fig:testCaseEKF}, the tracking error of rotor 159 is larger than the one of centralized EKF due to the lack of global information. But the difference is still small and within the acceptable range.

\subsection{Real PMU Data from Jiangsu Electric Power Company}
Jiangsu Electric Power Company, one of the largest provincial power company in China, has installed generation capacity of 100GW and peak load of 92GW. Over 160 PMUs, with thousands of measurement channels, have been installed in the Jiangsu system. These PMUs cover all 500kV substations, a majority of the 220kV substations, major power plants, and all renewable power plants. In this subsection, PMU data collected from real system oscillation events are used to validate the proposed oscillation monitoring algorithms.

\subsubsection{Case 1} 
\begin{figure}[ht]
\centering
\subfigure[Generator power measurements]{
\label{fig:JiangsuCase0Power}
  \includegraphics[width=.5\textwidth]{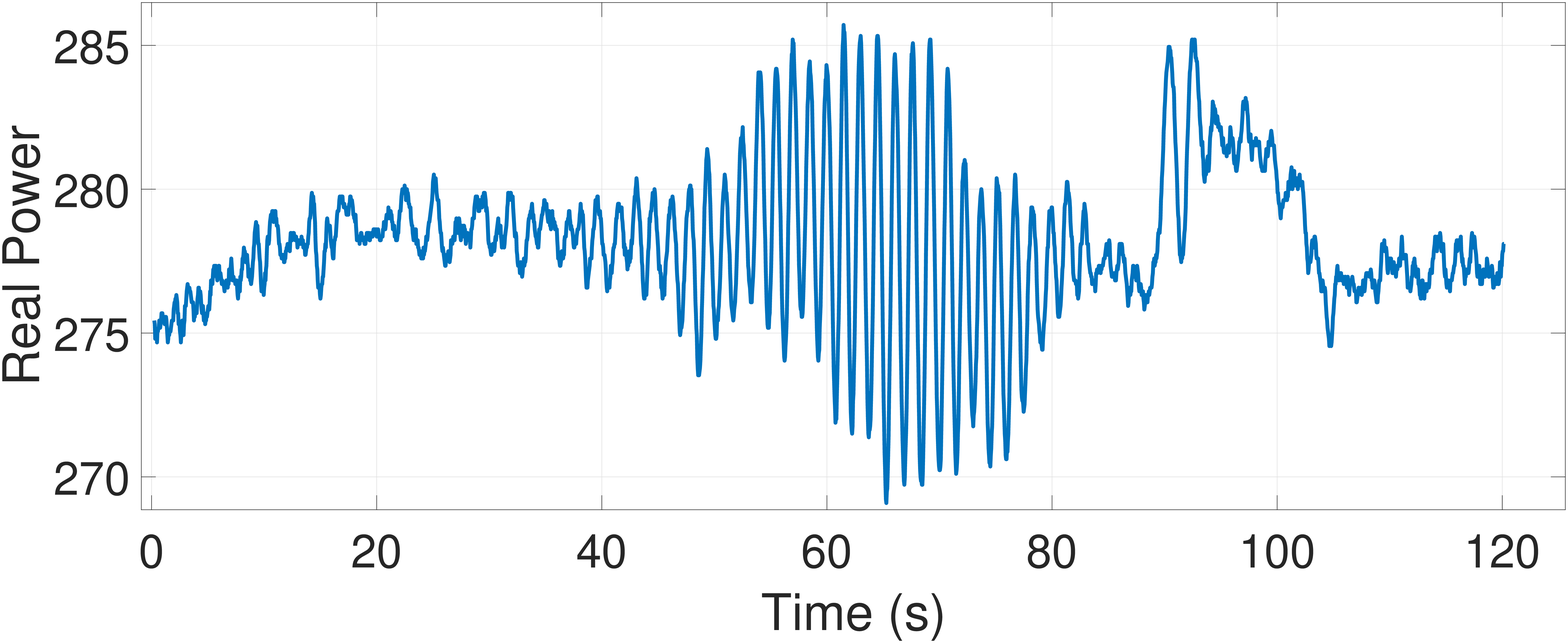}
  }
\subfigure[FFT results]{
\label{fig:JiangsuCase0FFT}
  \includegraphics[width=.5\textwidth]{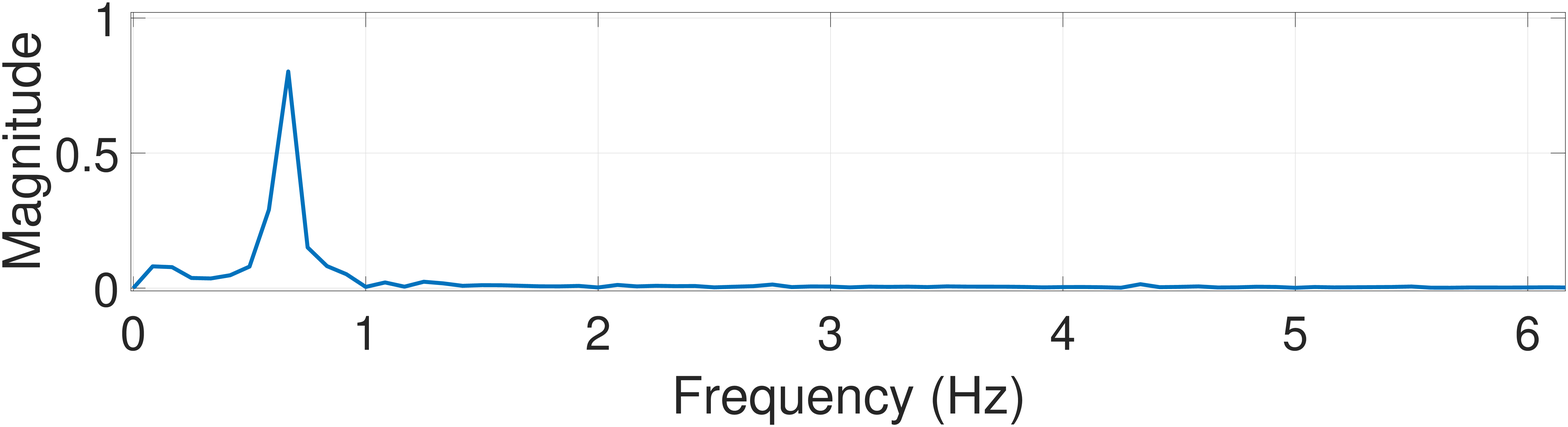}
  }
  \subfigure[Fitting of power]{
\label{fig:JiangsuCase0PowerFitting}
  \includegraphics[width=.5\textwidth]{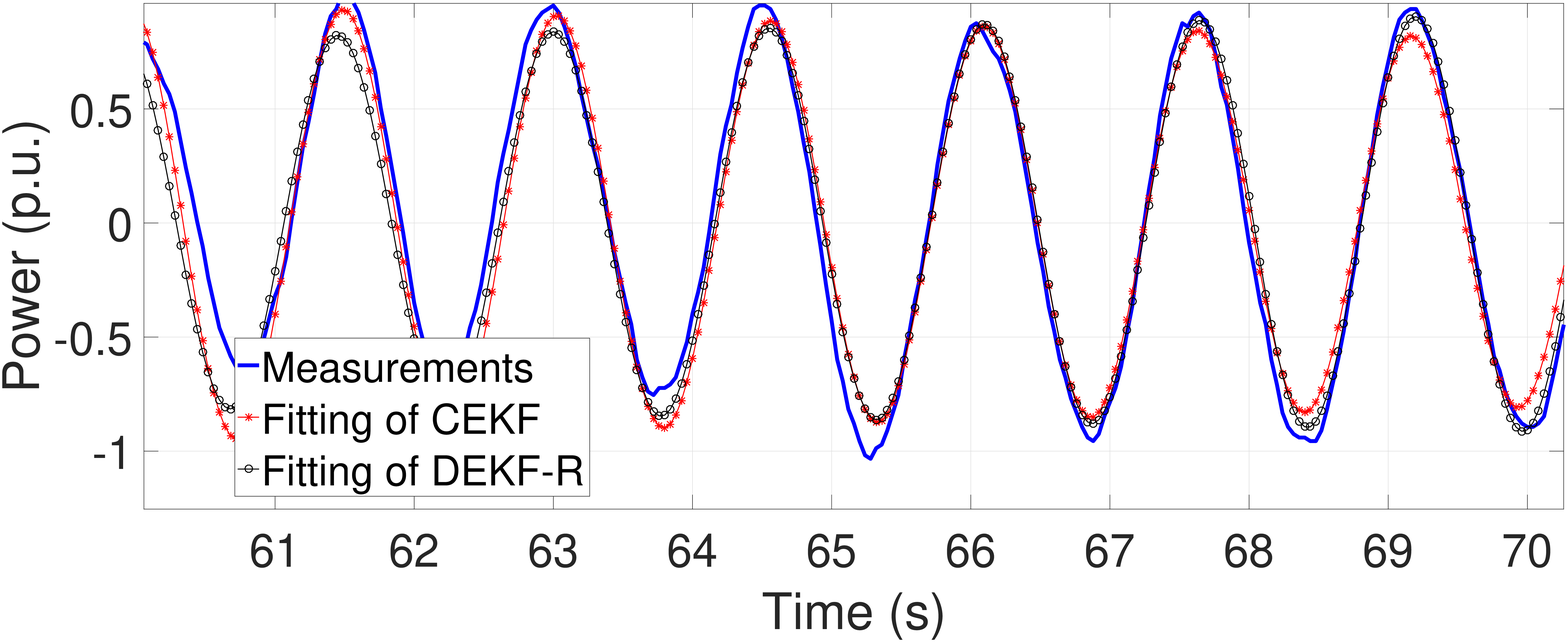}
  }
    \subfigure[Fitting of current magnitude]{
\label{fig:JiangsuCase0IMagFitting}
  \includegraphics[width=.5\textwidth]{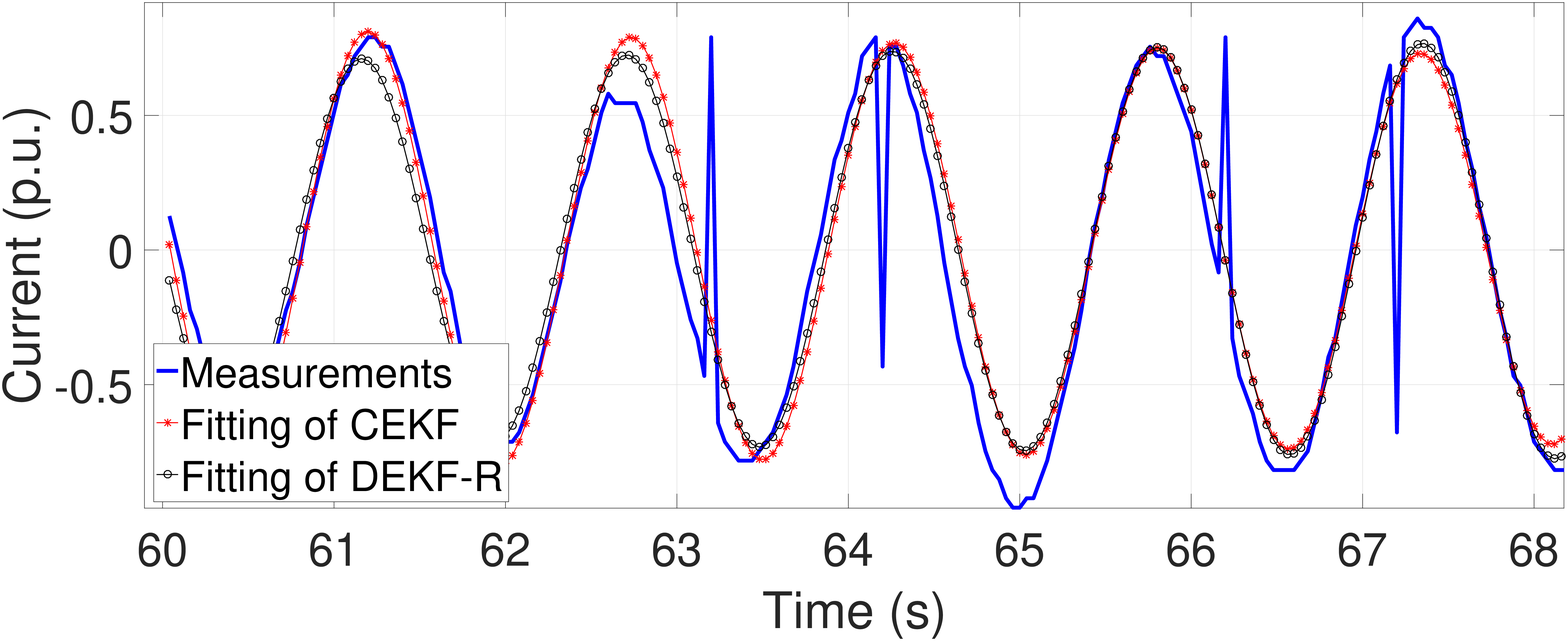}
  }
  \caption{Jiangsu data case 1: 1 oscillation mode}
\label{fig:JiangsuCase1}
\end{figure}

In this case, real and reactive power, voltage magnitude and angle, current magnitude and angle, and rotor speed and angle measurements from more than 160 PMUs are collected, at a reporting rate of $25$Hz. As shown in \ref{fig:JiangsuCase0Power}, an oscillation mode with frequency around $1$Hz lasts for 25 seconds. We remove the DC component in the data, carry out the normalization and feed the data from 60 to 80 seconds to FFT. As shown in Figure \ref{fig:JiangsuCase0FFT}, FFT results suggest  a dominant mode around $0.67$Hz. Using this result as initial estimates for Algorithm \ref{alg:CEKF} and \ref{alg:DEKF_reduce}, we estimate the damping factor and frequency in both centralized and distributed ways.

Estimation results are summarized in Table \ref{table:JiangsuCase1}. The estimated frequency of oscillation is around $0.65$ Hz with a damping factor of about $0.017$. The fitted curve of power and current magnitude of both centralized and distributed ways are presented in Figure \ref{fig:JiangsuCase0PowerFitting} and \ref{fig:JiangsuCase0IMagFitting}. Although errors exist in the current magnitude measurements, the fitting of both algorithm performs well.
\begin{table}
\caption{Jiangsu data case 1: 1 oscillation mode}
\centering
\begin{tabular}{|l|c|c|}
  \hline
 Estimate& $\omega$ & $\sigma$  \\
\hline
CEKF& $0.6513$ & $0.017$\\
\hline
DEKF-R& $0.6438$ & $0.0163$\\
\hline
\end{tabular}

\label{table:JiangsuCase1}
\end{table}

\subsubsection{Case 2}  
\begin{figure}[ht]
\centering
\subfigure[Generator voltage magnitude measurements]{
\label{fig:JiangsuCase3Vol}
  \includegraphics[width=.5\textwidth]{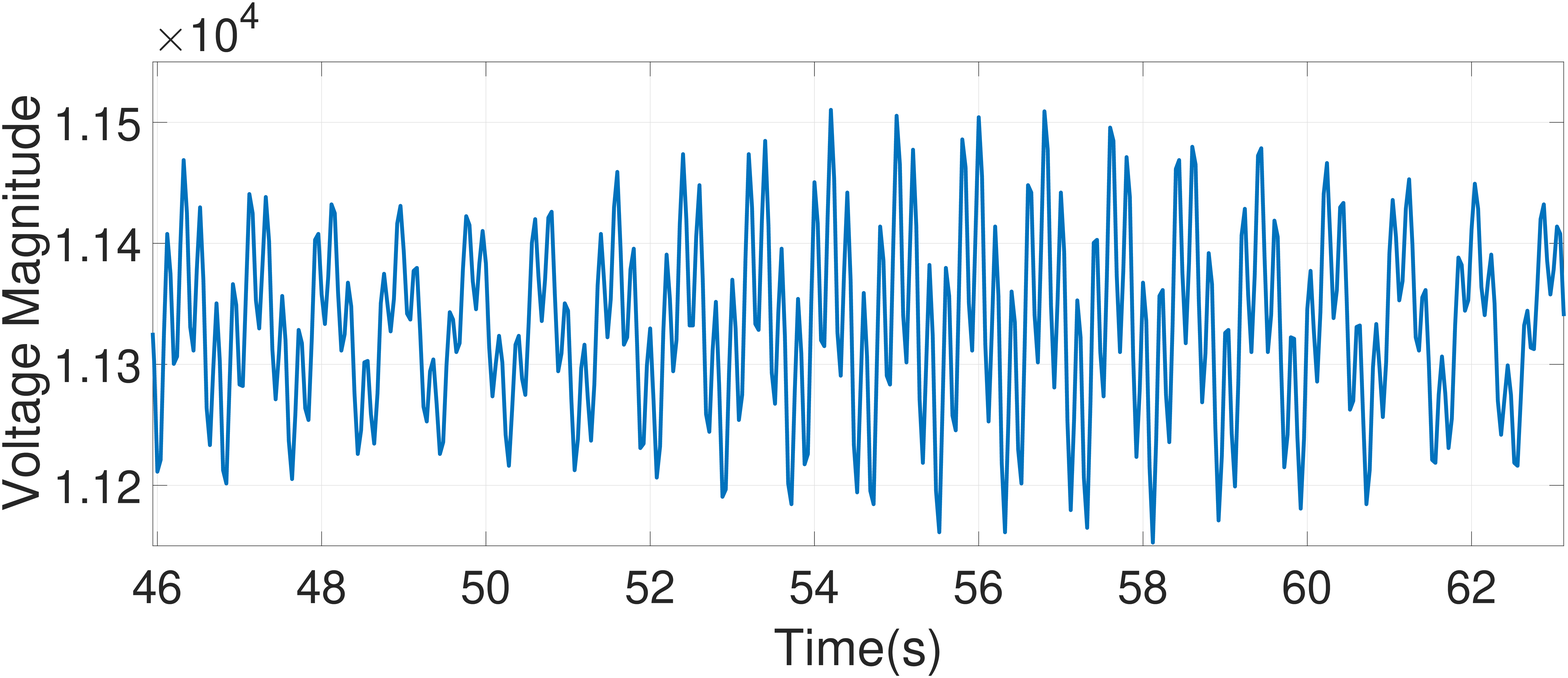}
  }
\subfigure[FFT results]{
\label{fig:JiangsuCase3FFT}
  \includegraphics[width=.5\textwidth]{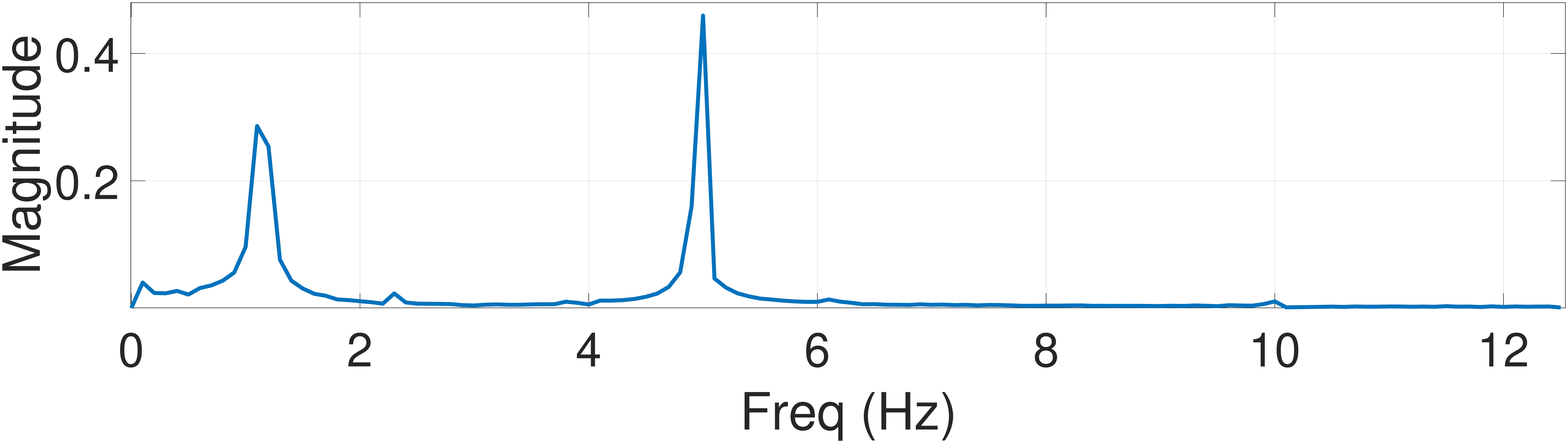}
  }
  \subfigure[Fitting of voltage magnitude]{
\label{fig:JiangsuCase3VolMagFitting}
  \includegraphics[width=.5\textwidth]{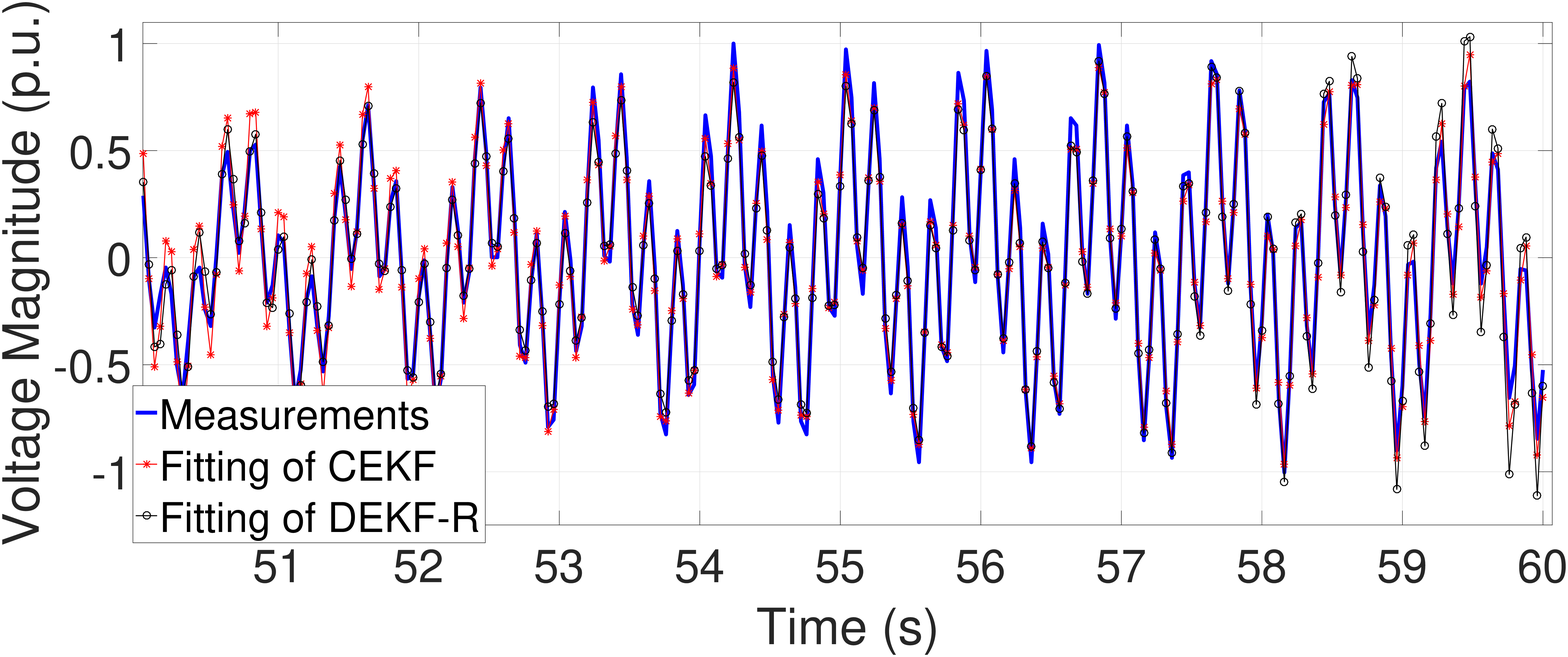}
  }
    \subfigure[Fitting of voltage phasor]{
\label{fig:JiangsuCase3VolAngFitting}
  \includegraphics[width=.5\textwidth]{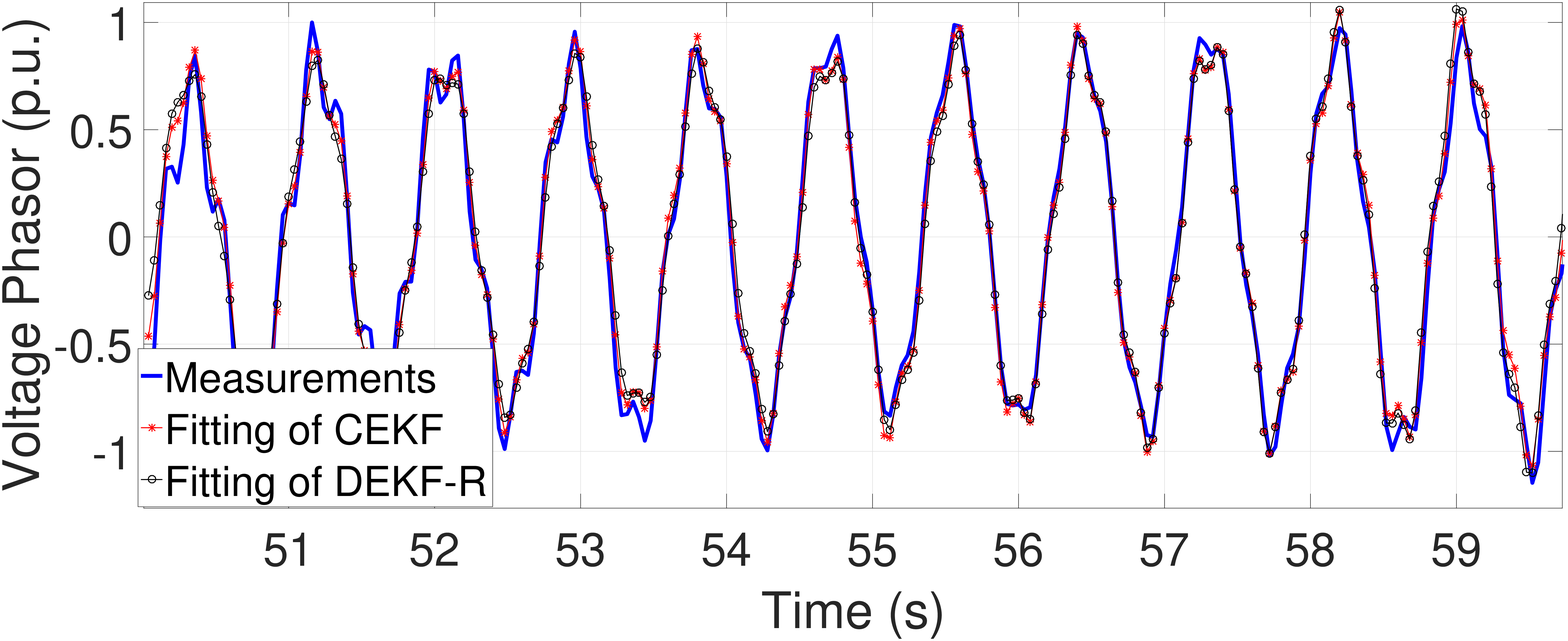}
  }
  \caption{Jiangsu data case 2: 2 oscillation modes}
\label{fig:JiangsuCase2}
\end{figure}

In this case, the voltage signal oscillates more than 60 seconds. Use a $10$ seconds window of data and repeat the same process. FFT suggests that there are two oscillation modes around 1Hz and 5Hz, respectively. Given this initial estimates, the proposed EKF algorithms are carried out, and estimation results can be found in Table \ref{table:JiangsuCase2} and Figure \ref{fig:JiangsuCase2}. The estimated frequencies are around $1.14$Hz and $4.97$Hz, with negative damping. Figure \ref{fig:JiangsuCase3VolMagFitting} and \ref{fig:JiangsuCase3VolAngFitting} show that the proposed methods perform well with multiple modes.

\begin{table}
\caption{Jiangsu Data case 2: 2 oscillation modes}
\centering
\begin{tabular}{|l|c|c|c|c|}
  \hline
 Estimate& $\omega$ & $\sigma$ & $\omega$ & $\sigma$ \\
\hline
CEKF& $1.1461$ & $-0.0215$& $4.9749$ & $-0.0222$\\
\hline
DEKF-R& $1.1450$ & $-0.0203$& $4.9494$ & $-0.0737$\\
\hline
\end{tabular}
\label{table:JiangsuCase2}
\end{table}

\section{Conclusion}\label{sec:VI}
Oscillation monitoring is essential in power systems to detect events and help system operators to identify the causes and locations of events. Conventional centralized algorithms put heavy burdens on data communications infrastructure and suffer from single point of failure and data privacy problem.  A novel distributed EKF based algorithm is proposed to estimate oscillation frequency and damping ratios directly.
 The fully distributed framework makes it possible to estimate at a fast reporting rate without information disclosure concerns. The effectiveness of the proposed algorithm is demonstrated using both simulated and real data.


{
\bibliographystyle{ieeetran}
\bibliography{Bibs/Journal,Bibs/Conf,Bibs/Book,Bibs/Misc}
}


\end{document}